\documentclass{article}
\usepackage{amsfonts}
\usepackage{amssymb}
\usepackage{amsbsy}
\usepackage{amsmath}
\usepackage{amsthm} 
\usepackage{comment}

\newtheorem{theorem}{Theorem}

\newtheorem{definition}[theorem]{Definition}

\newtheorem*{planX}{Plan for the proof of Theorem~\ref{the-main-theorem}}

\newcommand{\C}{\mathbf{C}}

\newcommand{\Q}{\mathbf{Q}}

\newcommand{\Z}{\mathbf{Z}}

\newcommand{\N}{\mathbf{N}}

\newcommand{\M}{\mathcal{M}}

\begin{document}
\title{On the slopes of the~$U_5$ operator acting on overconvergent modular forms}
\author{L. J. P. Kilford}
\maketitle
\begin{abstract}

We show that the slopes of the~$U_5$ operator acting on slopes of 5-adic overconvergent modular forms of weight~$k$ with primitive Dirichlet character~$\chi$ of conductor~25 are given by either
\[
\left\{\frac{1}{4}\cdot\lfloor\frac{8i}{5}\rfloor: i \in \N\right\}\text{ or }\left\{\frac{1}{4}\cdot\lfloor\frac{8i+4}{5}\rfloor: i \in \N\right\},
\]
depending on~$k$ and~$\chi$.
\end{abstract}

\section{Introduction}

We first define the slope of a (normalised) cuspidal eigenform.

\begin{definition}
Let~$f$ be a normalised cuspidal modular eigenform with $q$-expansion at~$\infty$ given by~$\sum_{n=1}^\infty a_n q^n$. The  \emph{slope} of~$f$ is defined to be the $5$-valuation of~$a_5$ viewed as an element of~$\C_5$; we normalise the $5$-valuation of~$5$ to be~$1$.
\end{definition}

As a consequence of the main result of this paper, we will prove the following theorem about classical modular forms.
\begin{theorem}
\label{classical-theorem}
The cyclotomic polynomial~$\Phi_{20}(x)$ factors over~$\Q_5$ into two factors, such that
\begin{eqnarray}
\label{2}
f_1 &\equiv& x^4 + 2x^3 + 4x^2 + 3x + 1 \mod 5,\\
\label{3}
f_2 &\equiv& x^4 + 3x^3 + 4x^2 + 2x + 1 \mod 5.
\end{eqnarray}

Let~$\chi$ be an odd primitive Dirichlet character of conductor~$25$ and let~$\tau$ be an odd primitive Dirichlet character of conductor~$5$.

Let~$k$ be a positive integer. We fix an embedding of the field of definition of~$\chi$ into~$\Q_5(\sqrt[4]{5},\sqrt{3})$.

Then the slopes of the~$U_5$ operator acting on~$S_k(\Gamma_0(25),\chi\tau^{k-1})$ are
\begin{eqnarray*}
&\left\{\frac{1}{4}\cdot\lfloor\frac{8i}{5}\rfloor: i \in \N\right\}&\text{ if }\chi(6)\text{ is a root of }f_1,\\
&\left\{\frac{1}{4}\cdot\lfloor\frac{8i+4}{5}\rfloor: i \in \N\right\}&\text{ if }\chi(6)\text{ is a root of }f_2.\\
\end{eqnarray*}
\end{theorem}

\section{Some previous work}

This paper uses methods introduced by Emerton in his PhD thesis~\cite{emerton-thesis}, which deals with the action of the~$U_2$ operator. It also uses methods developed by Smithline in his thesis~\cite{smithline}, which were then also used in the author's paper~\cite{kilford-2slopes} and in the paper of the author with Buzzard~\cite{buzzard-kilford}.

In~\cite{smithline-published}, the following theorem is proved about 3-adic modular forms:\begin{theorem}[Smithline~\cite{smithline-published}, Theorem~4.3]
We order the slopes of~$U_3$ by size, beginning with the smallest.

The sum of the first~$x$ nonzero slopes of the~$U_3$ operator acting on 3-adic overconvergent modular forms of weight~0 is at least~$3x(x-1)/2+2x$, and is exactly that if~$x$ is of the form~$(3^j-1)/2$ for some~$j$.
\end{theorem}
His thesis also shows that the sum of the first~$x$ slopes of the~$U_5$ operator acting on 5-adic overconvergent modular forms of weight~0 is at least~$x^2$.

In Buzzard-Kilford~\cite{buzzard-kilford}, the following theorem was proved about the 2-adic slopes of~$U_2$ acting on certain spaces of modular forms.
\begin{theorem}[Buzzard-Kilford~\cite{buzzard-kilford}, Theorem~$B$]
Let~$k$ be an integer and let~$\theta$ be a character of conductor~$2^n$ such that~$\theta(-1)=(-1)^k$.

If~$|5^k\cdot\theta(5)-1|_2>1/8$, then
the slopes of the overconvergent cuspidal modular forms of weight~$k$ and character~$\theta$ are $\{t,2t,3t,\ldots\}$, where $t=v(5^k\cdot\theta(5)-1)$, and
each slope occurs with multiplicity~1.
\end{theorem}

In Buzzard-Calegari~\cite{buzzard-calegari}, the following theorem is proved:
\begin{theorem}
The slopes of the~$U_2$ operator acting on 2-adic overconvergent modular forms of weight~0 are
\[
\left\{1 + 2v_2\left(\frac{(3n)!}{n!}\right): n \in \mathbf{N}\right\},
\]
where~$v_2$ is the normalised 2-adic valuation.
\end{theorem}

\section{Defining 5-adic overconvergent modular forms}

We now present the definition of the 5-adic overconvergent modular forms, first by defining overconvergent modular forms of weight~$0$, and then by deriving the definition for forms with weight and character.

This section follows Section~3 of~\cite{kilford-2slopes} in its layout and direction; more details on the specific steps can be found there.

Following Katz~\cite{katz}, section~2.1, we recall that, for~$C$ an elliptic curve over an $\mathbf{F}_5$-algebra~$R$, there is a mod~$5$ modular form~$A(C)$ called the \emph{Hasse invariant}, which has $q$-expansion over~$\mathbf{F}_5$ equal to~$1$.

We consider the Eisenstein series of weight~$4$ and tame level~1 defined over~$\mathbf{Z}$, with~$q$-expansion
\[
E_4(q):=1+240 \sum_{n = 1}^\infty \left( \sum_{0 < d | n} d^3 \right) \cdot q^n.
\]
We see that~$E_4$ is a lifting of~$A(C)$ to characteristic~0, as the reduction of~$E_4$ to characteristic~$5$ has the same $q$-expansion as~$A(C)$, and therefore~$E_4 \mod 5$ and~$A(C)$ are both modular forms of level~1 and weight~4 defined over~$\mathbf{F}_5$, with the same $q$-expansion. Note also that if~$C$ is an elliptic curve defined over~$Z_5$ then the valuation~$v_5(E_4(C))$ can be shown to be well-defined.

It is interesting to note that one can use the same Eisenstein series, $E_4$, in this part of the definition for 2-adic, 3-adic and 5-adic overconvergent modular forms (as a lifting of the $4^{th}$, $2^{nd}$ and~$1^{st}$ power of the Hasse invariant, respectively).

We now let~$m$ be a positive integer. Using arguments exactly similar to those in~\cite{kilford-2slopes}, we define the affinoid subdomain~$Z_0(5^m)$ of~$X_0(5^m)$ to be the connected component containing the cusp~$\infty$ of the set of points~$t=(C,P)$ in~$X_0(5^m)$ which have~$v_5(E_4(t))=0$. 

We now define strict affinoid neighbourhoods of~$Z_0(5^m)$.

\begin{definition}[Coleman~\cite{coleman}, Section~B2]
\label{connected-component}
We think of~$X_0(5^m)$ as a rigid space over~$\mathbf{Q}_5$, and we let~$t \in X_0(5^m)(\overline{\mathbf{Q}}_5)$ be a point, corresponding either to an elliptic curve defined over a finite extension of~$\Q_5$, or to a cusp.
Let~$w$ be a rational number, such that~$0 < w < 5^{2-m}/6$.

We define~$Z_0(5^m)(w)$ to be the connected component of the affinoid
\[
\left\{t \in X_0(5^m): \; v_5(E_4(t)) \le w\right\}
\]
which contains the cusp~$\infty$.
\end{definition}

Given this definition, we can now define 5-adic overconvergent modular forms.
\begin{definition}[Coleman,~\cite{coleman-overconvergent}, page~397]
Let~$w$ be a rational number, such that~$0 < w < 5^{2-m}/6$. Let~$\mathcal{O}$ be the structure sheaf of~$Z_0(5^m)(w)$.
We call sections of~$\mathcal{O}$ on~$Z_0(5^m)(w)$ 
\emph{$w$-overconvergent 5-adic modular forms of weight~$0$ and level~$\Gamma_0(5^m)$}.
If a section~$f$ of~$\mathcal{O}$ is a $w$-overconvergent modular form, then we say that~$f$ is an \emph{overconvergent $5$-adic modular form}.

Let~$K$ be a complete subfield of~$\mathbf{C}_5$, and define~$Z_0(5^m)(w)_{/K}$ to be the affinoid over~$K$ induced from~$Z_0(5^m)(w)$ by base change from~$\mathbf{Q}_5$. The space$$M_0(5^m,w;K):=\mathcal{O} (Z_0(5^m)(w)_{/K})$$ of $w$-overconvergent modular forms of weight~$0$ and level~$\Gamma_0(5^m)$ is a $K$-Banach space.

We now let~$\chi$ be a primitive Dirichlet character of conductor~$5^m$ and let~$k$ be an integer such that~$\chi(-1)=(-1)^k$. Let~$E^*_{k,\chi}$ be the normalised Eisenstein series of weight~$k$ and character~$\chi$ with nonzero constant term.

The space of $w$-overconvergent 5-adic modular forms of weight~$k$ and character~$\theta$ is given by
\[
\M_{k,\theta}(5^m,w;K):=E^*_{k,\theta} \cdot \M_0(5^m,w;K).
\]
This is a Banach space over~$K$.
\end{definition}

There are Hecke operators~$U_5$ and~$T_p$ (where~$p\nmid 5$) acting on the space of modular forms~$\M_{k,\theta}(5^m,w;K)$; these are defined on the $q$-expansions of the overconvergent modular forms in exactly the same way as they are defined on the $q$-expansions of classical modular forms. One defines~$T_n$ for~$n$ a natural number in the usual way.

Using results of Coleman, we have the following theorem about the independence of the characteristic power series of~$U_5$ acting on~$\M_{k,\theta}(5^m,w;K)$:
\begin{theorem}[Coleman~\cite{coleman}, Theorem~B3.2]
Let~$w$ be a real number such that~$0 < w < \min(5^{2-m}/6,1/6)$, let~$k$ be an integer and let~$\theta$ be a character such that~$\theta(-1)=(-1)^k$.

The characteristic polynomial of~$U_5$ acting on $w$-overconvergent 5-adic modular forms of weight~$k$ and character~$\theta$ is independent of the choice of~$w$.
\end{theorem}

We will now rewrite the definition of~$Z_0(25)(w)$ in terms of a carefully chosen modular function of level~$25$, in order to prove the following theorem:
\begin{theorem}
Let~$w_0=1/12$. The space of $w_0$-overconvergent modular forms of weight~0 and level~$25$, with coefficients in~$\Q_5(\sqrt[4]{5})$, is a Tate algebra in one variable over~$\Q_5(\sqrt[4]{5})$.
\end{theorem}
\begin{proof}
We have given a valuation on the points~$t$ of the rigid space~$X_0(5^m)$, based on the lifting of the Hasse invariant by the Eisenstein series~$E_4$. We recall that the modular $j$-invariant is defined to be~$j:=E_4^3/\Delta$. Therefore, we see that, if the elliptic curve corresponding to~$t$ has good reduction, then~$\Delta(t)$ has valuation~0, and therefore that 
\[
v_5(t)=v_5(E_4(t))=\frac{1}{3} v_5((E_4)(t)^3)=\frac{1}{3} v_5(j(t)).
\]
We now recall that the modular curve~$X_0(25)$ has genus~0. This means that there is a modular function~$t_{25}$ which is a uniformiser on~$X_0(25)$: 
\[
t_{25}:=\frac{\eta(q)}{\eta(q^{25})},
\]
where~$\eta$ is the Dedekind~$\eta$-function.
We could write~$t_{25}$ as a rational function in~$j$ directly, but as the resulting rational function is very complicated, we will instead also work with the uniformiser~$t_5$ of~$X_0(5)$, defined as
\[
t_5:=\left(\frac{\eta(q)}{\eta(q^5)}\right)^6
\]
By explicit calculation, one can verify the following identities of modular functions:
\begin{eqnarray}
\label{rat-fns}
j=\frac{(t_5^2+250t_5+3125)^3}{t_5^5}\text{ and }t_5=\frac{t_{25}^5}{t_{25}^4+5t_{25}^3+15t_{25}^2+25t_{25}+25}.
\end{eqnarray}
We note also that
\[
j(\infty)=t_5(\infty)=t_{25}(\infty)=\infty;
\]
this follows because the $q$-expansion of all of these functions begins~$q^{-1}+\cdots$.

Because~$t_{5}(\infty)=t_{25}(\infty)=\infty$, the connected components of~$Z_0(5)$ and~$Z_0(25)$ which contain the cusp~$\infty$ are of the form~$v_5(t_5) < D_1$ and~$v_5(t_{25}) < D_2$, for some rational numbers~$D_1$ and~$D_2$.

By considering the Newton polygons of the numerators and denominators of the rational functions in~\eqref{rat-fns}, we see that if~$v_5(t_{25}) < 1/2$, then~$v_5(t_{25})=v_5(t_5)=v_5(j)$. This means that we have shown that
\[
Z_0(25)(w)=\left\{x \in X_0(25):\;v_5(t_{25}(x))\le 3w\right\}\text{, for }0 < w < 1/6.
\]
Now, we choose~$w=1/12$, and therefore we obtain
\[
Z_0(25)(1/12)=\left\{x \in X_0(25):\;v_5(t_{25}(x))\le 1/4\right\}.
\]
Let us define~$W:=\sqrt[4]{5}/t_{25}$. We can rewrite the definition of~$Z_0(25)(1/12)$ again in terms of~$W$ to get
\[
Z_0(25)(1/12)=\left\{x \in X_0(25):\;v_5(W(x))\ge 0\right\}.
\]

Finally, we recall that the rigid functions on the closed disc over~$\mathbf{Q}_5$ with centre~0 and radius~1 are defined to be power series of the form
\[
\sum_{n \in \mathbf{N}}a_n z^n\;:\;a_n \in \mathbf{Q}_5,\;a_n \rightarrow 0.
\]
Therefore, the $1/12$-overconvergent modular forms of level~$\Gamma_0(25)$ and weight~0 are 
\[
\mathbf{Q}_5(\sqrt[4]{5})\langle W \rangle,
\]
which is what we wanted to show.
\end{proof}

We have written down this space of overconvergent modular forms as an explicit Banach space. This means that we can write down its \emph{Banach basis}: the set~$\left\{W,W^2,W^3,\ldots\right\}$ forms a Banach basis for the overconvergent modular forms of weight~$0$ and level~$\Gamma_0(25)$. This Banach basis is composed of weight~$0$ modular functions --- we want to be able to consider the action of the~$U_5$ operator on overconvergent modular forms with non-zero weight~$k$ and character~$\chi$ (here, as elsewhere in this note, $\chi$ has conductor~25 and~$\theta(-1)=(-1)^k$). Using an observation from the work of Coleman~\cite{coleman}, we will be able to move between weight~0 and weight~$k$ and character~$\chi$ via multiplication by a suitable quotient of modular forms.

Let~$F$ be an overconvergent modular form of weight~$k$ and character~$\theta$ which has nonzero constant term, and let~$z$ be an overconvergent modular function of weight~$0$. In particular, we note that~$F$ may have negative weight. From the discussion in Coleman~\cite[page~450]{coleman} we see that the pullback~$\tilde{U}_5$ of the~$U_5$ operator acting on overconvergent modular forms of weight-character~$(k,\theta)$ to weight~$0$ is~$1/F\cdot U_5(z\cdot F)$. 

Now by equation~3.3 of~\cite{coleman-old} we have that~$U_5(z\cdot V(F))=U_5(z)\cdot F$. We therefore consider the modular form~$H=V(G)$, and substitute~$H$ for~$F$ in the formula we have just derived for~$U_2(z\cdot V(F))$, to obtain:
\[
\tilde{U}_5(z\cdot V(G))=\frac{1}{V(G)}\cdot U_5(z\cdot V(G))=\frac{G}{V(G)}\cdot U_5(z).
\]
We can also use this line of reasoning to see that
\[
1/F\cdot U_5(z\cdot F)=U_5\left({z\cdot\frac{F}{V(F)}}\right).
\]
This allows us to now define the (twisted)~$U_5$ operator.
\begin{definition}[The twisted~$U_5$ operator]
\label{twisted-U-definition}
Let~$k$ be an integer and let~$\chi$ be an odd character of conductor~25. Let~$\tau$ be a character of conductor~5 such that~$\chi(-1)=(-1)^k$.

The twisted~$U_5$ operator acting on forms of weight-character~$(1+kt,\chi\cdot\tau^{k})$ is defined to be the following operator:
\begin{equation}
\label{horrible-u5}
U_5\left({W^i \cdot \frac{E^*_{1,\chi}}{V(E^*_{1,\chi})}}\right)\cdot\left(\frac{E^*_{k,\tau}}{V(E^*_{k,\tau})}\right)^t.
\end{equation}
\end{definition}

We can now consider the action of this twisted~$U_5$ operator on these spaces of overconvergent modular forms.
\begin{definition}[The matrix of the twisted~$U_5$ operator]
\label{twisted-matrix-definition}
Let~$k$ be an integer and let~$\chi$ be an odd character of conductor~25. Let~$\tau$ be a character of conductor~5 such that~$\chi(-1)=(-1)^k$.

Let~$M=(m_{i,j})$ be the infinite compact matrix of the twisted~$U_5$ operator acting on overconvergent modular forms of weight-character~$(1+kt,\chi\cdot\tau^{k})$, 
where~$m_{i,j}$ is defined to be the coefficient of~$W^i$ in the $W$-expansion of the operator defined in equation~\eqref{horrible-u5}.
\end{definition}
Here we will make the observation that the entries of our matrix~$M$ are functions of~$\chi$, $\tau$, and~$t$.

We know that~$U_5$ is a compact operator, so 
we can show 
that the trace, determinant and characteristic power series of~$M$ are all well-defined. We will use a theorem of Serre to prove our theorem on the slopes of~$U_5$ acting on~$M$.
\begin{theorem}[Serre~\cite{serre}, Proposition~7]
\label{serre-proposition}

\begin{enumerate}
\item Let~$M_n$ be an~$n \times n$ matrix defined over a finite extension of~$\Q_2$. Let~$\mathop{det}(1-tM_n)=\sum_{i= 0}^n c_i t^i$. Let~$M_m$ be the matrix formed by the first~$m$ rows and columns of~$M_n$.

Let~$s(i)$ be the formula for the~$i^{th}$ slope; in our specific case, this will mean that either
\[
s(i)=\frac{1}{4}\cdot\lfloor\frac{8i}{5}\rfloor{\rm \text{ or }}s(i)=\frac{1}{4}\cdot\lfloor\frac{8i+4}{5}\rfloor.
\]

Assume that there exists a constant~$r \in \Q^\times$ such that 
\begin{enumerate}
\item For all positive integers~$m$ such that~$1 \le m \le n$, the valuation of~$\mathop{det}(M_m)$ is~$r \cdot \sum_{i=1}^m s(i)$.

\item The valuation of elements in column~$j$ is at least~$r \cdot s(i)$.
\end{enumerate}

Then we have that, for all positive integers~$m$ such that~$1 \le m \le n$, $v_2(c_m)=r \cdot \sum_{i=1}^m s(i)$.

\item Let~$M_\infty$ be a compact infinite matrix (that is, the matrix of a compact operator). If~$M_m$ is a series of finite matrices which tend to ~$M_\infty$, then the finite characteristic power series~$\mathop{det}(1-tM_m)$ converge coefficientwise to~$\mathop{det}(1-tM_\infty)$, as~$m \rightarrow \infty$.

\end{enumerate}
\end{theorem}

We now quote a result of Coleman that tells us that overconvergent modular forms of small slope are in fact classical modular forms:
\begin{theorem}[Coleman~\cite{coleman-overconvergent}, Theorem~1.1]
\label{coleman-ocgt}
Let~$k$ be a non-negative integer and let~$p$ be a prime. Every $p$-adic overconvergent modular eigenform of weight~$k$ with slope strictly less than~$k-1$ is a classical modular form.
\end{theorem}

We now state the main theorem of this paper, which tells us exactly what the slopes of the~$U_5$ operator acting on slopes of modular forms of level~25.
\begin{theorem}
\label{the-main-theorem}
We recall and use the notation of Theorem~\ref{classical-theorem}.


Let~$\chi$ be an odd primitive Dirichlet character of conductor~$25$ and let~$\tau$ be an odd primitive Dirichlet character of conductor~$5$.

Let~$k$ be a positive integer. We fix an embedding of the field of definition of~$\chi$ into~$\Q_5(\sqrt[4]{5})$, and recall the notation of~$f_1$ and~$f_2$ from Theorem~\ref{classical-theorem}.

The slopes of overconvergent modular forms of weight~$k$ and character~$\chi\tau^{k-1}$ are given by
\begin{eqnarray*}
&\left\{\frac{1}{4}\cdot\lfloor\frac{8i}{5}\rfloor: i \in \N\right\}&\text{ if }\chi(6)\text{ is a root of }f_1,\\
&\left\{\frac{1}{4}\cdot\lfloor\frac{8i+4}{5}\rfloor: i \in \N\right\}&\text{ if }\chi(6)\text{ is a root of } f_2.\\
\end{eqnarray*}
\end{theorem}

We can prove Theorem~\ref{classical-theorem}, assuming Theorem~\ref{the-main-theorem}, by recalling the following theorem from Cohen-Oesterl\'e:
\begin{theorem}[Cohen-Oesterl\'e~\cite{cohen-oesterle}, Th\'eor\`eme~1]
Let~$\chi$ be a primitive Dirichlet character of conductor~25 and let~$k$ be a positive integer greater than~1. The following formula holds:
\[
d(k,\chi):=\dim S_k(\Gamma_0(25),\chi)=\frac{5k-7}{2}+\varepsilon\cdot (\chi(8)+\chi(17)),
\]
where~$\varepsilon$ is~0 for odd~$k$, $-1/4$ if~$k \equiv 2 \mod 4$, and~$1/4$ if~$k \equiv 0 \mod 4$.
\end{theorem}

\begin{proof}[Proof of Theorem~\ref{classical-theorem}]
The classical theorem will follow, because when we substitute~$d(k,\chi)$ into the formula~$s(i)$ for the~$i^{th}$ slope, we see that the maximum value of~$s(d(k,\chi))$ is~$k-1$. We now apply either Theorem~\ref{coleman-ocgt} or an argument of Buzzard shows that slopes greater than~$k-1$ cannot be classical (see~\cite{buzzard-kilford}, the proof of the Corollary to Theorem~B, which references the proof of Theorem~4.6.17(1) of~\cite{miyake}), so therefore as there are at most~$k-1$ slopes which are smaller than or equal to~$k-1$, we see that all of these small slopes are the slopes of classical eigenforms.
\end{proof}
\section{Observations}

There are some interesting new features which appear when~$p=5$ that do not appear when~$p=2$ or~$p=3$.

Firstly, there is a computational issue. In previous work (such as~\cite{buzzard-kilford} or~\cite{kilford-2slopes}), the computations in {\sc magma}~\cite{magma} could be carried out either over the rational numbers or over the field~$\Q_p$. However, for~$p=5$ the calculation must be carried out over~$\Q_5(\sqrt[4]{5})$.

Secondly, and more importantly, there are now two different possibilities for the slopes, which depend on which character is chosen. These two possibilities correspond to the two factors of the cyclotomic polynomial~$\Phi_{20}(x)$ over~$\Q_5$. This is a departure from the situation in~\cite{buzzard-kilford} and~\cite{kilford-2slopes}, where the slopes are independent of choice of character.

Thirdly, the slopes are no longer in one arithmetic progression, as they are in the previously studied cases. Instead, there are five arithmetic progressions which interlace together; these all have a common difference between terms (which is~2). (One could, of course, view the arithmetic progression~$1,2,3,\ldots$ as being made up of the two arithmetic progressions~$1,3,5,\ldots$ and~$2,4,6\ldots$ but this point of view is only reasonable after one has considered the action of~$U_5$ on forms of level~25, where the slopes form several arithmetic progressions).

Fourthly, part of the complexity in the calculations in~\cite{kilford-2slopes} was the fact that (in the notation of that paper) the modular function~$U_2(z^{2i+1})$ was identically zero. This meant that the first ``matrix of the~$U_2$ operator'' that was defined had identically zero determinant, which meant that some algebra had to be done to get a matrix to which Theorem~\ref{serre-proposition} could be applied to. In the current note, this does not happen because the matrix~$M$ of the~$U_5$ operator does not have any identically zero columns. 

Finally, the strategy of Section~\ref{technical-horrors} was chosen because the modular functions involved were unusually simple, thus making the calculations more tractable (the corresponding functions for (say) weight~2 were much less pleasant).
 
\section{The technical part; proof of Theorem~\ref{the-main-theorem}}
\label{technical-horrors}

As the actual proof of Theorem~\ref{the-main-theorem} is somewhat technical, we will first outline a plan to show how the proof works.

\begin{planX}
\emph{In this section, we will show that we can apply Theorem~\ref{serre-proposition}, which will prove Theorem~\ref{the-main-theorem}. 
First we fix an arbitrary positive integer~$n$, an integer~$k$ and a primitive Dirichlet character~$\theta$ of conductor~$25$ such that~$\theta(-1)=-1$.} 

\emph{We will begin with the matrix~$M_{n}$; the matrix formed by the first~$n$ rows and~$n$ columns of~$M$, the matrix of the twisted~$U_5$ operator acting on forms of weight-character~$(1,\theta)$ defined in Definition~\ref{twisted-U-definition}. The proof will then proceed in the following way:}

\begin{enumerate}

\item Define the matrix~$D(\beta(i))$ to be the diagonal matrix with~$\beta(j)$ in the~$j$th row and the~$j$th column. We define the matrix~$O_n:=D(\sqrt{5}^{-j})\cdot O_n\cdot D(\sqrt{5}^j)$.

\item We then show that the valuation of elements in the~$j$th column of~$O_n$ are~$s(j)$; this verifies condition~$(b)$ of Theorem~\ref{serre-proposition}, with~$r=s(j)$.

\item We finally show that~$O_n$ has determinant of valuation~$\sum_{i=1}^n s(i)$, by considering the matrix~$P_n:=D(5^{-s(j)})\cdot O_n$. By showing that~$P_n$ has determinant of valuation~0, it can be seen that the valuation of the determinant of~$O_n$ is the valuation of the determinant of~$D(5^{s(j)})$, which is~$\sum_{i=1}^n s(i)$. This will verify condition~$(a)$ of Theorem~\ref{serre-proposition}, with determinant of valuation~$\sum_{i=1}^n s(i)$.

\item Finally, we will show that, after multiplication by the multiplier (as defined in Definition~\ref{twisted-matrix-definition})
\[
\left(\frac{E^*_{1,\tau}}{V(E^*_{1,\tau})}\right)^t,
\]
the matrix of the twisted~$U_5$ operator acting on forms of weight~$1+t$ and character~$\theta\cdot\tau^{t-1}$ still satisfies properties~$(a)$ and~$(b)$ of Theorem~\ref{serre-proposition}.
\end{enumerate}

\emph{At each step of this plan, we must show that the characteristic polynomial of the new matrix defined is the same as that of~$M_n$. In the last step, we will show that~$P_n$ has unit determinant by reducing it modulo a prime ideal above~$5$ and showing that this reduction has determinant~$1$. This means that we must prove that~$P_n$ has coefficients which are integers in~$\Q_5(\sqrt[4]{5})$.}

\end{planX}

\begin{proof}[Proof of Theorem~\ref{the-main-theorem}]

In this section, we will use the modular function~$T$ instead of~$W$; we define~$T$ as follows:
\[
T:=\frac{1}{t_{25}}.
\] This will make it easier to perform the calculations.

By computation, it can be shown that the first five columns of the matrix~$M_n$ (in weight~1) are polynomials in~$T$ of degree~$5i$; we will just give the valuations of the coefficients of these, as the actual coefficients are elements of~$\Q_5(\sqrt[4]{5})$ and thus take up a lot of space.

If we have chosen~$\chi(6)$ to be a root of~\eqref{2}, then these are the valuations:
\begin{eqnarray*}
U_5(T\cdot E^*_{1,\chi}):&\left[ \frac{1}{2}, \frac{3}{2}, \frac{9}{4}, \frac{13}{4}, 4 \right] \\
U_5(T^2\cdot E^*_{1,\chi}):&\left[ \frac{1}{2}, 1, 2, \frac{11}{4}, 4, \frac{9}{2}, \frac{11}{2}, \frac{25}{4}, \frac{29}{4}, 8 \right]\\
U_5(T^3\cdot E^*_{1,\chi}):&\left[ \frac{1}{4}, 1, \frac{5}{4}, \frac{9}{4}, 3, 5, 5, 6, \frac{27}{4}, 8, \frac{17}{2}, \frac{19}{2}, \frac{41}{4}, \frac{45}{4}, 12 \right]\\
U_5(T^4\cdot E^*_{1,\chi}):& \left[ \frac{1}{4}, \frac{3}{4}, \frac{5}{4}, \frac{7}{4}, 3, \frac{7}{2}, \frac{9}{2}, \frac{21}{4}, \frac{25}{4}, 7, \frac{17}{2}, 9, 10, \frac{43}{4}, 12, \frac{25}{2}, \frac{27}{2}, \frac{57}{4}, \frac{61}{4}, 16 \right]\\
U_5(T^5\cdot E^*_{1,\chi}):& [ 0, 1, 1, 2, 2, \frac{7}{2}, 4, 5, \frac{23}{4}, 7, \frac{15}{2}, \frac{17}{2}, \frac{19}{2}, \frac{21}{2}, 11, \frac{25}{2}, 13, 14, \frac{59}{4}, 16, \frac{33}{2}, \frac{35}{2}, \frac{73}{4}, \frac{77}{4}, 20 ]
\end{eqnarray*}
If, on the other hand, we have chosen~$\chi(6)$ to be a root of~\eqref{3}, then these are the valuations:
\begin{eqnarray*}
U_5(T\cdot E^*_{1,\chi}):&\left[\frac{1}{4},\frac{5}{4},2,3,4\right] \\
U_5(T^2\cdot E^*_{1,\chi}):&\left[ \frac{1}{4}, \frac{3}{4}, \frac{7}{4}, \frac{5}{2}, 4, \frac{17}{4}, \frac{21}{4}, 6, 7, 8 \right] \\
U_5(T^3\cdot E^*_{1,\chi}):&\left[0, \frac{3}{4}, 1, 2, 3, \frac{17}{4}, \frac{19}{4}, \frac{23}{4}, \frac{13}{2}, 8, \frac{33}{4}, \frac{37}{4}, 10, 11, 12
\right] \\
U_5(T^4\cdot E^*_{1,\chi}):&\left[ 0, \frac{1}{2}, 1, \frac{3}{2}, 3, \frac{13}{4}, \frac{17}{4}, 5, 6, 7, \frac{17}{2}, \frac{35}{4}, \frac{39}{4}, \frac{21}{2}, 12, \frac{49}{4}, \frac{53}{4}, 14, 15, 16\right] \\
U_5(T^5\cdot E^*_{1,\chi}):&\left[0, 1, 1, 2, 2, \frac{7}{2}, \frac{15}{4}, \frac{19}{4}, \frac{11}{2}, 7, \frac{29}{4}, \frac{33}{4}, \frac{37}{4}, \frac{41}{4}, 11, \frac{49}{4}, \frac{51}{4}, \frac{55}{4}, \frac{29}{2}, 16, \frac{65}{4}, \frac{69}{4}, 18, 19, 20 
\right]. 
\end{eqnarray*}

The valuations of the~$T$-coefficients of~$U_5(T^5)$ are independent of the choice of~$\chi(6)$ (because~$T$ has $q$-coefficients in~$\Q_5$) and are as follows:\[
U_5(T^5):\;\left[ 0, \frac{1}{4}, \frac{1}{4}, \frac{1}{2}, \frac{1}{2}, \frac{5}{4}, \frac{3}{2}, \frac{3}{2}, \frac{7}{4}, \frac{7}{4}, \frac{5}{2}, \frac{5}{2}, \frac{5}{2}, \frac{11}{4}, \frac{11}{4}, \frac{7}{2}, 4, \frac{15}{4}, 4, 4, \frac{9}{2}, \frac{19}{4}, \frac{19}{4}, 5, \frac{19}{4} \right].
\]

We also note that the following identity of modular functions holds:
\begin{equation}
\label{weight-1-mult}
\frac{E^*_{1,\tau}}{V(E^*_{1,\tau})}=1-\frac{5(T+(2+2I)T^2)}{1+(2+I)T+(2+I)T^2}.\end{equation}

\subsection{Checking the valuations of the elements in the~$j^{th}$ column of~$O_n$} 
We will show that the valuations of elements in the~$jth$ column of the matrix~$O_n$ have valuation at least~$s(j)$ by writing the operator which determines the~$(5a+b)^{th}$ column, $U_5(T^{5a+b} \cdot E^*_{1,\chi}/V(E^*_{1,\chi}))$ 
in terms of the operators~$U_5(T^b\cdot E^*_{1,\chi}/V(E^*_{1,\chi}))$ and~$U_5(T
^5)$.

\begin{definition}
We will write
\[
\overline{U_5(T\cdot E^*_{1,\chi}/V(E^*_{1,\chi}))}
\]
to mean that we have changed the basis of the matrix~$M_n$ of the twisted~$U_5$ operator by conjugating it with the matrices~$D(5^{s(j)})$ and~$D(5^{-s(j)})$.

We define this matrix to be~$O_n$, as defined in the Plan given at the beginning of this section.
\end{definition}

Let~$\pi$ be a fixed fourth root of~5. By consulting the tables of valuations above, we see that (ignoring unit factors in the coefficients) the following congruences hold:
\begin{eqnarray}
\label{u5t1}
\overline{U_5(T\cdot E^*_{1,\chi}/V(E^*_{1,\chi}))} &\equiv& \pi^{4s(1)} T\mod \pi^{4s(1)+1}\\
\overline{U_5(T^2\cdot E^*_{1,\chi}/V(E^*_{1,\chi}))} &\equiv&\pi^{4s(2)} (T+T^2) \mod 
\pi^{4s(2)+1}\\
\overline{U_5(T^3\cdot E^*_{1,\chi}/V(E^*_{1,\chi}))} &\equiv&\pi^{4s(3)}(T+T^3) \mod \pi^{4s(3)+1}\\
\overline{U_5(T^4\cdot E^*_{1,\chi}/V(E^*_{1,\chi}))} &\equiv& \pi^{4s(4)}(T+T^2+T^3+T^
4)\mod \pi^{4s(4)+1}\\
\overline{U_5(T^5\cdot E^*_{1,\chi}/V(E^*_{1,\chi}))} &\equiv& \pi^{4s(5)}(T+T^3+T^5) \mod \pi^{4s(5)+1}\\
\label{u5t5}\overline{U_5(T^5)} &\equiv& \pi^8(T+T^3+T^5)\mod \pi^9.
\end{eqnarray}

Now we see that
\begin{equation}
\label{taking-the-T5-out}
\overline{U_5\left(T^{5a+b}\cdot \frac{E^*_{1,\chi}}{V(E^*_{1,\chi})}\right)} \equiv \overline{U_5\left(T^{b}\cdot \frac{E^*_{1,\chi}}{V(E^*_{1,\chi})}\right)}\cdot \overline{U_5(T^5)}^a;
\end{equation}
this is because~$T^5$ is congruent to~$V(T)$ modulo~5, and from the definition of the~$U_5$ and~$V$ operators, we see that~$\overline{U_5(T^{5a+b}\cdot X)}$ (where~$X$ is a modular function) is congruent to~$T^a \cdot \overline{U_5(T^{b}\cdot X)} $ modulo~5. Finally, one can check explicitly that the $q$-expansions of~$T$ and~$\overline{U_5(T^5)}$ are congruent modulo~5. 

Now we can show that the valuations of the~$T$-coefficients of~$\overline{U(T^{5a+b}\cdot E^*_{1,\chi}/V^*_{1,\chi})}$ are at least~$s(5a+b)$, as is required for condition~$(b)$ of Theorem~\ref{serre-proposition}; from the argument above we see that
\[
\overline{U_5(T^{5a+b}\cdot X)} \equiv \overline{U_5(T^5)}^a \cdot \overline{U_5(T^{b}\cdot X)} \mod 5.
\]
This means that the following equality holds:
\begin{equation}
\label{whydoesthiswork}
\overline{U_5(T^{5a+b}\cdot X)} = \overline{U_5(T^5)}^a \cdot \overline{U_5(T^{b}\cdot X)} + 5^\varepsilon f(T),
\end{equation}
for some function~$f(T)$ with integral coefficients in some extension of~$\Z_5$.

Now we know from the discussion above that the valuation of~$T$-coefficients of~$\overline{U_5(T^5)}$ is at least~2, and that the valuation of~$T$-coefficients of~$\overline{U_5(T^{b}\cdot X)}$ is at least~$s(b)$, so therefore the valuation of the $T$-coefficients of~$\overline{U_5(T^{5a+b}\cdot X)}$ is at least~$s(5a+b)$, as is required.


\subsection{Defining~$P_n$ and showing that it has unit determinant}
We now postmultiply the matrix~$O_n$ by~$D(5^{-s(i)})$ and define this product to be~$P_n$. We will now show that~$P_n$ has determinant of valuation~0, and therefore that the valuation of the determinant of~$O_n$ is~$\sum_{i=1}^n s(i)$.

We now reduce the entries of the matrix~$P_n$ modulo a prime above~5; we call this matrix~$P^\prime_n$. Now, to show that~$P_n$ has unit determinant it will suffice to show that the columns of~$P^\prime_n$ are linearly independent.

From the congruences shown in~\eqref{u5t1}--\eqref{u5t5} and the argument at~\eqref{taking-the-T5-out} and~\eqref{whydoesthiswork}, we see that the elements of the~$(5a+b)^{th}$ column of~$P^\prime_n$ are given by the coefficients~$T$ in the~$T$-expansion of
\[
\overline{U_5\left(T^{b}\cdot \frac{E^*_{1,\chi}}{V(E^*_{1,\chi})}\right)}\cdot\overline{U_5(T)}^a,
\]
and that the highest coefficient of this that does not vanish after reduction modulo the prime ideal above~5 is exactly~$T^{5a+b}$.

We can now see immediately that the matrix~$P^\prime_n$ has determinant a unit, because there are units on the diagonal, and no elements below the diagonal. Therefore the determinant of~$P_n$ is a 5-adic unit, as required.

This means that the determinant of~$O_n$ and also the determinant of~$M_n$ both have valuation~$\sum_{i=1}^n s(i)$. This means that~$M_n$ satisfies condition~$(a)$ of Theorem~\ref{serre-proposition} and therefore that we can apply this theorem to the matrix~$M$ to show that the slopes of the~$U_5$ operator are given by~$s(i)$.

\subsection{Generalising all this to other weights}
We note that this part of the proof has shown that the matrices~$M_n$ of the twisted~$U_5$ operator acting on overconvergent modular forms of weight~1 have determinants with 5-valuation~$\sum_{i=1}^m s(i)$, and that the valuations of elements in the~$j^{th}$ column are at least~$s(i)$. We will now prove that the matrix of the twisted~$U_5$ operator acting on weights of the form~$1+t$, where~$t$ is an integer, also satisfies these two properties; this will be enough to prove Theorem~\ref{the-main-theorem}.

It should be noted that the strange-seeming conditions involving~$\chi$ and~$\tau$ in Theorem~\ref{classical-theorem} and Theorem~\ref{the-main-theorem} arise because we first proved the theorems for weight~1 and will then use the multiplier~$(E^*_{1,\tau}/V(E^*_{1,\tau}))^{k-1}$ to get to weight~$k$.

We now check that the multiplier given in~\eqref{weight-1-mult} has large enough valuations that after the change-of-basis, its $T$-coefficients still have non-negative valuation.

Under the embedding we have chosen of~$\Q_5(\chi)$ into the extension field of~$\Q_5$, both~$3-I$ and~$2+I$ have normalised valuation~1. This means that the valuation of the coefficient of~$T^2$ in the multiplier has valuation~1, so therefore under the change of basis by conjugation by the two diagonal matrices, we see that the valuation of the coefficients of the multiplier will be non-negative. This means that we are still able to reduce modulo a prime ideal above~5, so our analysis will carry through.

We now define matrices~$O_n$, $P_n$ and~$P^\prime_n$ in the same way as above, but now using the operator
\[
U_5\left(T^{5a+b}\cdot\frac{E^*_{1,\chi}}{V(E^*_{1,\chi})}\right)\cdot \left(\frac{E^*_{1,\tau}}{V(E^*_{1,\tau})}\right)^t.
\]
Because we have checked that the valuations of the weight~1 and level~5 multiplier~\eqref{weight-1-mult} are compatible with those of the original twisted~$U_5$ operator in weight~1, we merely need to check that the columns of the matrix~$P^\prime_n$ are still linearly independent.

We notice that, after the change of basis and reduction modulo the prime ideal, the weight~1 multiplier is of the form~$1+\cdots$; in other words, it is a unit. So the columns of~$P^\prime_n$ in weight~$1+t$ are linearly independent, because the columns of the matrix~$P^\prime_n$ in weight~1 are linearly independent, and we have multiplied each of these columns by a unit.

This means that~$P_n$ has unit determinant, and therefore that~$O_n$ and~$M_n$ have determinant of valuation~$\sum_{i=1}^n s(i)$. Therefore we can apply Theorem~\ref{serre-proposition} and hence we have proved Theorem~\ref{the-main-theorem}.
\end{proof}
\section{Acknowledgements}
A preliminary version of this article was written during my stay at the Max-Planck-Institut in Bonn, Germany in the summer of 2004. The author would like to thank the Institut for its hospitality.

The results of this paper were mentioned in a short note circulated at the Banff International Research Station workshop ``$p$-adic variation of motives'', which was held in December 2003. The author would also like to thank the BIRS for its hospitality.

Some of the computer calculations were executed by the computer algebra package {\sc Magma}~\cite{magma} running on the machine \textbf{crackpipe}, which was bought by Kevin Buzzard with a grant from the Central Research Fund of the University of London. I would like to thank the CRF for their support. Other computer calculations took place on William Stein's machine~\textbf{Meccah} at Harvard University; I would like to thank him for the use of his computer.

I would like to thank Kevin Buzzard, Edray Goins and Ken McMurdy for helpful conversations.


\begin{thebibliography}{10}

\bibitem{magma}
W.~Bosma, J.~Cannon, and C.~Playoust.
\newblock The {M}agma algebra system {I}: The user language.
\newblock {\em J. Symb. Comp.}, 24(3--4):235--265, 1997.
\newblock \emph{http://magma.maths.usyd.edu.au}.

\bibitem{buzzard-calegari}
Kevin Buzzard and Frank Calegari.
\newblock Slopes of overconvergent 2-adic modular forms.
\newblock {\em Compos. Math.}, 141(3):591--604, 2005.

\bibitem{buzzard-kilford}
Kevin Buzzard and L.~J.~P. Kilford.
\newblock The 2-adic eigencurve at the boundary of weight space.
\newblock {\em Compos. Math.}, 141(3):605--619, 2005.

\bibitem{cohen-oesterle}
H.~Cohen and J.~Oesterl\'e.
\newblock Dimensions des espaces de formes modulaires.
\newblock {\em Lecture Notes in Mathematics}, 627:69--78, 1977.

\bibitem{coleman-overconvergent}
R.~Coleman.
\newblock Classical and overconvergent modular forms of higher level.
\newblock {\em J. Th\'eor. Nombres Bordeaux}, 9(2):395--403, 1997.

\bibitem{coleman}
R.~Coleman.
\newblock $p$-adic {B}anach spaces and families of modular forms.
\newblock {\em Inv. Math}, 127:417--479, 1997.

\bibitem{coleman-old}
Robert~F. Coleman.
\newblock Classical and overconvergent modular forms.
\newblock {\em Invent. Math.}, 124(1-3):215--241, 1996.

\bibitem{emerton-thesis}
M.~Emerton.
\newblock {\em 2-adic Modular Forms of minimal slope}.
\newblock PhD thesis, Harvard University, 1998.

\bibitem{katz}
N.~Katz.
\newblock $p$-adic properties of modular forms and modular curves.
\newblock {\em Lecture Notes in Mathematics}, 350:69--190, 1973.

\bibitem{kilford-2slopes}
L.~J.~P. Kilford.
\newblock Slopes of 2-adic overconvergent modular forms with small level.
\newblock {\em Math. Res. Lett.}, 11(5-6):723--739, 2004.

\bibitem{miyake}
T.~Miyake.
\newblock {\em Modular Forms}.
\newblock Springer, 1989.

\bibitem{serre}
J.-P. Serre.
\newblock Endomorphismes completements continues des espaces de {B}anach
  $p$-adique.
\newblock {\em Publ. Math. IHES}, 12:69--85, 1962.

\bibitem{smithline}
L.~Smithline.
\newblock {\em Exploring slopes of $p$-adic modular forms}.
\newblock PhD thesis, University of California at Berkeley, 2000.

\bibitem{smithline-published}
Lawren Smithline.
\newblock Compact operators with rational generation.
\newblock In {\em Number theory}, volume~36 of {\em CRM Proc. Lecture Notes},
  pages 287--294. Amer. Math. Soc., Providence, RI, 2004.

\end{thebibliography}

\end{document}